\documentclass[suppldata]{interact}

\usepackage{latexsym}
\usepackage{calc}
\usepackage{amssymb}
\usepackage{xcolor}
\usepackage[center]{caption}
\usepackage{hyperref}
\hypersetup{colorlinks=true, allcolors=blue,
pdfinfo={Title={Gradient Methods with Memory},
Author={Yurii Nesterov and Mihai I. Florea},
Subject={Gradient Methods for Convex Minimization},
Keywords={convex optimization, gradient methods, relative smoothness, rate of convergence, piece-wise linear model}}}

\usepackage[numbers,sort&compress]{natbib}
\bibpunct[, ]{[}{]}{,}{n}{,}{,}

\def\Def{\stackrel{\mathrm{def}}{=}}

\def\inter{{\rm int \,}}

\def\dom{{\rm dom \,}}

\def\R{\mathbb{R}}
\def\E{\mathbb{E}}

\newcommand{\refLE}[1]{\overset{\eqref{#1}}\leq}
\newcommand{\refEQ}[1]{\overset{\eqref{#1}}=}
\newcommand{\refGE}[1]{\overset{\eqref{#1}}\geq}

\newcommand{\refGT}[1]{\overset{\eqref{#1}}>}

\newtheorem{theorem}{Theorem}
\newtheorem{lemma}{Lemma}

\newtheorem{assumption}{Assumption}

\newtheorem{remark}{Remark}
\newcommand{\half}{\mbox{${1 \over 2}$}}

\def\la{\langle}
\def\ra{\rangle}

\begin{document}

\title{Gradient Methods with Memory\footnote{This is an Accepted Manuscript of an article published by Taylor \& Francis in Optimization Methods and Software on 13 Jan 2021, available at \url{https://www.tandfonline.com/doi/10.1080/10556788.2020.1858831}.}}

\addtocounter{footnote}{-1}

\author{
\name{Yurii Nesterov\textsuperscript{a}
and Mihai I. Florea\textsuperscript{b}
\thanks{CONTACT Mihai I. Florea. E-mail: \href{mailto:mihai.florea@uclouvain.be}{mihai.florea@uclouvain.be}}
}
\affil{\textsuperscript{a} Center for Operations Research and
Econometrics (CORE), Catholic University of Louvain (UCL), Belgium;
\textsuperscript{b} Department of Mathematical Engineering (INMA),
UCL, Belgium}
}

\maketitle

\begin{abstract}
In this paper, we consider gradient methods for
minimizing smooth convex functions, which employ the
information obtained at the previous iterations in order
to accelerate the convergence towards the optimal
solution. This information is used in the form of
a piece-wise linear model of the objective function, which
provides us with much better prediction abilities as
compared with the standard linear model. To the best of
our knowledge, this approach was never really applied in
Convex Minimization to differentiable functions in view of
the high complexity of the corresponding auxiliary
problems. However, we show that all necessary computations
can be done very efficiently. Consequently, we get new
optimization methods, which are better than the usual
Gradient Methods both in the number of oracle calls and
in the computational time. Our theoretical conclusions are
confirmed by preliminary computational experiments.
\end{abstract}


\section{Introduction}

\subsection{Motivation} First-order gradient methods for
minimizing smooth convex functions generate a sequence
of test points based on the information obtained from the
oracle: the function values and the gradients. Most
methods either use the information from the last test
point or accumulate it in the form of an aggregated
linear function (see, for example, Chapter 2 in
\cite{ref_007}). This approach is very different from the
technique used in Nonsmooth Optimization, where the
piece-wise linear model of the objective is a standard and
powerful tool. It is enough to mention the Bundle Method,
the Level Method, cutting plane schemes, etc. The reason for
this situation is quite clear. The presence of piece-wise
linear models in the auxiliary problems, which we need to
solve at each iteration of the method, usually
significantly increases the complexity of the corresponding
computations. This is acceptable in Nonsmooth
Optimization, which has the reputation of a difficult field.
By contrast, Smooth Optimization admits very simple
and elegant schemes, with a very small computational cost of
each iteration, preventing us from introducing there such
a heavy machinery.

After the preparation of this manuscript, we became aware of
a highly specialized attempt in \cite{ref_002},
using quadratic lower bounds instead of linear ones.
Although the results presented there seem promising, the study limits
itself to studying smooth unconstrained problems with strongly convex
objectives. The study also states that an extension of the results
to a wider context is a difficult open problem.

The main goal of this paper is the demonstration that the
above situation is not as clear as it looks like. We will
show that the Gradient Method,\footnote{By Gradient Method we denote
the extended scheme described in \cite{ref_001}, which encompasses
Gradient Descent and the Proximal Gradient Method.} equipped with a
piece-wise linear model of the objective function, has much more
chances to accelerate on particular optimization
problems. At the same time, it appears that the
corresponding auxiliary optimization problems can be
easily solved by an appropriate version of the Frank-Wolfe
algorithm. All our claims are supported by a complexity
analysis. In the end, we present preliminary computational
results, which show that very often the new schemes have
much better computational time.

\subsection{Contents} In Section \ref{sc-Max}, we analyse the
Gradient Method with Memory as applied to the composite
form of smooth convex optimization problems
\cite{ref_005}. In order to measure the level of
smoothness of our objective function, we introduce the
{\em relative smoothness} condition~\cite{ref_001,ref_004}, based
on an arbitrary strictly convex distance function. The main
novelty here is the piece-wise linear model of the
objective function, formed around the current test point.
We analyse the corresponding auxiliary optimization
problem and propose a condition on its approximate
solution which does not destroy the rate of convergence of
the algorithm. In Section \ref{sc-Approx}, we analyse the
complexity of solving the auxiliary optimization problem
using the Frank-Wolfe algorithm. More precisely, we consider
the anti-dual\footnote{The dual of the problem with the objective
multiplied by $-1$.} of the auxiliary problem.

In this section, we restrict ourselves to
strongly convex distance functions. We show that our
auxiliary problem can be easily solved by the Frank-Wolfe
method. Its complexity is proportional to the maximal
squared norm of the gradient in the current model of the
objective divided by the desired accuracy.

In Section \ref{sc-Euclid}, we specify our complexity
results for the Euclidean setup, when all distances are
measured by a Euclidean norm. We show that, for some
strategies of updating the piece-wise linear model, the
complexity of the auxiliary computations is very low.

Finally, in Section \ref{sc-Num} we present preliminary
computational results. We compare the usual Gradient
Method with two gradient methods with memory, which use
different strategies for updating the piece-wise model of
the objective function. Our conclusion is that the new
schemes are always better, both in terms of the calls of
oracle and in the total computational time.

\subsection{Notation and generalities} In what follows, denote
by $\E$ a finite-dimensional real vector space and by
$\E^{*}$ its dual space, the space of linear functions on
$\E$. The value of function $s \in \E^{*}$ at point $x \in
\E$ is denoted by $\langle s, x \rangle$. Let us fix some
arbitrary (possibly non-Euclidean) norm $\| \cdot \|$ on
the space $\E$ and define the dual norm $\| \cdot \|_{*}$
on $\E^{*}$ in the standard way:
$$
\|s\|_{*} \Def \sup_{h \in \E} \{ \la s, h \ra : \|h\|
\leq 1 \}.
$$

Let us choose a simple closed convex
prox-function $d(\cdot)$, which is differentiable at the
interior of its domain.\footnote{Recall that a
function is closed if its epigraph is a closed set.}
This function must be strictly convex:
\begin{equation}\label{eq-SConv}
d(y) > d(x) + \la \nabla d(x), y - x \ra, \quad x \in
\inter(\dom d), \; y \in \dom d, \; x \neq y.
\end{equation}
Using this function, we can define the {\em Bregman
distance} between two points $x$ and $y$:
\begin{equation}\label{def-Breg}
\beta_d(x,y) = d(y) - d(x) - \la \nabla d(x), y - x \ra,
\quad x \in \inter(\dom d), \; y \in \dom d.
\end{equation}
Clearly, $\beta_d (x,y) \refGT{eq-SConv} 0$ for $x \neq y$
and $\beta_d(x,x) = 0$.

We will use Bregman distances for measuring the level of
{\em relative smoothness} of convex functions (see
\cite{ref_004}). Namely, for a differentiable closed convex
function $f$ with open $\dom f \subseteq \dom d$ we define two
constants, $L_d(f) \geq \mu_d(f) \geq 0$, such that
\begin{equation}\label{def-RSmooth}
\begin{array}{rcll}
f(y) - f(x) - \la \nabla f(x), y - x \ra & \geq & \mu_d(f)
\beta_d(x,y),\\
& & & \quad x,y \in \dom f,\\
f(y) - f(x) - \la \nabla f(x), y - x \ra & \leq & L_d(f)
\beta_d(x,y),
\end{array}
\end{equation}
See \cite{ref_001} and \cite{ref_004} for definitions,
motivations, and examples.

\section{Gradient method with memory}\label{sc-Max}

In this paper, we are solving the following composite
minimization problem:
\begin{equation} \label{prob-Main}
\min_{x \in \dom \psi} \Bigl\{ F(x) \equiv f(x) + \psi(x)
\Bigr\},
\end{equation}
where function $f$ satisfies the relative smoothness
condition (\ref{def-RSmooth}), possibly with $\mu_d(f) =
0$. The function $\psi: \E \to \R \cup \{ +\infty \}$ is a
proper closed convex function with $\dom \psi \subseteq \dom f$ and
$\inter (\dom \psi)$ non-empty. The function $\psi$ is simple
(in the sense of satisfying Assumptions \ref{ass-Phi} and
\ref{ass-Comp}, stated in the sequel) but it does not have to be
differentiable or even continuous. For instance, $\psi$ can incorporate
the indicator function of the feasible set.
We assume that a solution $x^{*} \in \dom \psi$ of problem
(\ref{prob-Main}) does exist, denoting $F^{*} = F(x^{*})$.

The simplest method for solving the problem
(\ref{prob-Main}) is the usual Gradient Method:
\begin{equation}\label{met-GM}
\begin{array}{l}
\mbox{Choose $x_0 \in \inter (\dom \psi)$. For $k \geq 0$, iterate:} \\
\\
x_{k+1} = \arg\min\limits_{y \in \dom \psi} \left\{ f(x_k)
+ \la \nabla f(x_k), y - x_k \ra + \psi(y) + L
\beta_d(x_k,y) \right\}.
\end{array}
\end{equation}
The constant $L$ in this method has to be big enough in
order to ensure
$$
\begin{array}{rcl}
f(x_{k+1}) & \leq & f(x_k) + \la \nabla f(x_k),
x_{k+1}-x_k \ra + L \beta_d(x_k, x_{k+1}).
\end{array}
$$
In view of (\ref{def-RSmooth}), this is definitely true
for $L \geq L_d(f)$. However, we are interested in
choosing $L$ as small as possible since this would
significantly increase the rate of convergence of the
scheme.

Method (\ref{met-GM}) is based on the simplest {\em linear
model} of function $f(\cdot)$ around the point $x_k$. In
our paper, we suggest to replace it by a piece-wise linear
model, defined by the information collected at other test
points.

Namely, for each $k \geq 0$ define a {\em discrete set}
${\cal Z}_k$ of $m_k$ feasible points ($m_k \geq 1$):
$$
{\cal Z}_k = \{ z_i \in \dom \psi, \; i = 1, \dots, m_k
\}.
$$
Then we can use a more sophisticated model of the smooth
part of the objective function,
\begin{equation}\label{eq-Mod}
f(y) \geq \ell_k(y) \Def \max\limits_{z_i \in {\cal Z}_k}
\{ f(z_i) + \la \nabla f(z_i), y - z_i \ra \}, \quad y \in
\dom \psi.
\end{equation}
This model is always better than the initial linear model
provided that
\begin{equation}\label{eq-Inc}
x_k \in {\cal Z}_k.
\end{equation}
In what follows, we always assume that this condition is
satisfied.

Thus, we come to the following natural generalization of
the method (\ref{met-GM}), which we call the {\em Gradient
Method with Memory} (GMM):
\begin{equation}\label{met-GMM}
\begin{array}{|l|}
\hline \\
\mbox{{\bf Choose} $x_0 \in \inter (\dom \psi)$.}\\
\\
\mbox{\bf For $k \geq 0$, iterate:} \\
\\
x_{k+1} = \arg\min\limits_{y \in \dom \psi} \left\{
\ell_k(y) + \psi(y) + L \beta_d(x_k,y) \right\}.\\
\\
\hline
\end{array}
\end{equation}
\begin{remark}\label{rm-L}
Note that for any $x \in \dom f$ we have
$$
\begin{array}{c}
f(x_k) + \la \nabla f(x_k), x-x_k \ra + L \beta_d(x_k, x)
\; \refLE{eq-Inc} \; \ell_k(x) + L
\beta_d(x_k, x) \\
\\
\refLE{eq-Mod} \; f(x) + L \beta_d(x_k, x) \;
\refLE{def-RSmooth} \; f(x_k) + \la \nabla f(x_k), x-x_k
\ra + (L+L_d(f)) \beta_d(x_k, x).
\end{array}
$$
Therefore, we can count on a better convergence of method
(\ref{met-GMM}) only if we will be able to choose the
parameter $L$ significantly smaller than $L_d(f)$.
\end{remark}

At each iteration of method (\ref{met-GMM}), we need to
solve a non-trivial auxiliary minimization problem.
Therefore, the practical efficiency of this method
crucially depends on the complexity of this computation.
In what follows, we suggest to solve this problem
approximately using a special method for its dual problem.

Let us start by presenting the corresponding technique.
For the sake of notation, we omit the index of the
iteration. Thus, our auxiliary problem is as follows:
$$
\min\limits_{y \in \dom \psi} \max\limits_{\lambda \in
\Delta_m} \left\{ \sum\limits_{i=1}^m \lambda^{(i)} [f_i +
\la g_i, y - z_i \ra] + \psi(y) + L \beta_d(\bar x, y )
\right\},
$$
where $f_i = f(z_i)$, $g_i = \nabla f(z_i)$, $i = 1,
\dots, m$, and $\Delta_m$ is the standard simplex in
$\R^m$. Introducing now the vector $f_* \in \R^m$ with
coordinates
\begin{equation}\label{def-f*}
f^{(i)}_* = \la g_i, z_i \ra - f_i, \quad i = 1, \dots, m,
\end{equation}
we get the following representation of our problem:
\begin{equation}\label{prob-Aux}
\min\limits_{y \in \dom \psi} \max\limits_{\lambda \in
\Delta_m} \left\{ \la \lambda, G^T y - f_* \ra + \psi(y) +
L \beta_d(\bar x, y ) \right\},
\end{equation}
where $G = (g_1, \dots, g_m) \in \E^* \times \R^m$. Note
that the pay-off function in this saddle point problem can
be written as follows:
$$
\begin{array}{rl}
& \la \lambda, G^T y - f_* \ra + \psi(y) + L \beta_d(\bar
x, y
)\\
\\
= & \la \lambda, G^T y - f_* \ra + \psi(y) + L [d(y) -
d(\bar x) - \la \nabla d(\bar x), y - \bar x \ra ]\\
\\
= & L d(y) - \la L \nabla d(\bar x) - G \lambda, y \ra +
\psi(y) - \la \lambda, f_* \ra + L [ \la \nabla d(\bar x),
\bar x \ra - d(\bar x)].
\end{array}
$$
Hence, we need to introduce the following dual function:
\begin{equation}\label{def-Phi}
\Phi_L(s) = \max\limits_{y \in \dom \psi} \left\{ \la s, y
\ra - L d(y) - \psi(y) \right\}, \quad s \in \E^*.
\end{equation}
Our main joint assumption on functions $d(\cdot)$ and
$\psi(\cdot)$ is as follows.
\begin{assumption}\label{ass-Phi}
For any $L>0$, function $\Phi_L(s)$ is defined at any $s
\in \E^*$.
\end{assumption}
This can be ensured, for example, by the strong convexity of
function $d(\cdot)$, or by the boundedness of $\dom \psi$,
or in many other ways.

Since the objective function in the definition
(\ref{def-Phi}) is strictly concave, its solution
$$
y^*_L(s) = \arg\max\limits_{y \in \dom \psi} \left\{ \la
s, y \ra - L d(y) - \psi(y) \right\}
$$
is uniquely defined for any $s \in \E^*$. Moreover,
function $\Phi_L(\cdot)$ is differentiable and
\begin{equation}\label{eq-GPhi}
\nabla \Phi_L(s) = y^*_L(s), \quad s \in \E^*.
\end{equation}

Now we can write down the problem anti-dual to
(\ref{prob-Aux})
\begin{equation}\label{prob-ADual}
\xi^*_L \; \Def \; \min\limits_{\lambda \in \Delta_m}
\left\{ \xi_L(\lambda) \Def \Phi_L(L\nabla d(\bar x) - G
\lambda) + \la \lambda, f_* \ra + \alpha \right\},
\end{equation}
where $\alpha = L [ d(\bar x) - \la \nabla d(\bar x), \bar
x \ra]$. This is a convex optimization problem with
a differentiable objective function. Our second main
assumption is as follows.
\begin{assumption}\label{ass-Comp}
Function $\Phi_L(\cdot)$ in the problem (\ref{prob-ADual})
is easily computable.
\end{assumption}

We will discuss the reasonable strategies for finding an
approximate solution to problem~(\ref{prob-ADual}) in
Section \ref{sc-Approx}. At this moment, it is enough to
assume that we are able to compute a point $\bar \lambda =
\bar \lambda (\bar x, {\cal Z}, L)$ such that
\begin{equation}\label{eq-OptL}
\begin{array}{rcl}
\la \bar \lambda - \lambda, \nabla \xi_L(\bar \lambda)
\ra& \leq & \delta, \quad \lambda \in \Delta_m,
\end{array}
\end{equation}
where $\delta \geq 0$ is some tolerance parameter.
Clearly, if $\delta = 0$, then $\bar \lambda$ is the
optimal solution to the problem (\ref{prob-ADual}). Note
that condition (\ref{eq-OptL}) ensures also a small
functional gap:
\begin{equation}\label{eq-XiGap}
\begin{array}{rcl}
\xi_L(\bar \lambda) - \xi_L^* & = & \max\limits_{\lambda
\in \Delta_m} [ \xi_L(\bar \lambda) - \xi_L(\lambda)] \;
\leq \; \max\limits_{\lambda \in \Delta_m} \la \bar
\lambda - \lambda, \nabla \xi_L(\bar \lambda) \ra \;
\refLE{eq-OptL} \; \delta.
\end{array}
\end{equation}

Condition (\ref{eq-OptL}) immediately leads to the
following result.
\begin{lemma}\label{lm-Tol}
Let $\bar \lambda \in \Delta_m$ satisfy condition
(\ref{eq-OptL}). Then for $\bar s = L \nabla d(\bar x) - G
\bar \lambda$ we have
\begin{equation}\label{eq-Tol}
\begin{array}{rcl}
\sum\limits_{i=1}^m \bar \lambda^{(i)} [ f_i + \la g_i,
y^*_L(\bar s) - z_i \ra] & \geq & \max\limits_{1 \leq i
\leq m}[ f_i + \la g_i, y^*_L(\bar s) - z_i \ra] - \delta.
\end{array}
\end{equation}
\end{lemma}
\begin{proof}
Indeed, $\nabla \xi_L(\bar \lambda) = f_* - G^T y^*_L(\bar
s)$. Thus, inequality (\ref{eq-OptL}) can be rewritten as
follows:
$$
\begin{array}{rcl}
\la \bar \lambda, G^T y^*_L(\bar s) - f_* \ra & \geq & \la
\lambda, G^T y^*_L(\bar s) - f_* \ra - \delta, \quad
\lambda \in \Delta_m.
\end{array}
$$
It remains to note that
$$
(G^T y^*_L(\bar s) - f_*)^{(i)} = f_i + \la g_i,
y^*_L(\bar s) - z_i \ra, \quad i = 1, \dots, m.$$
\end{proof}

Now we are able to analyse one iteration of the inexact
version of method (\ref{met-GMM}).
\begin{equation}\label{it-IGMM}
\begin{array}{|rl|}
\hline & \\
\mbox{\bf Input:} & \mbox{Point $\bar x \in \inter (\dom
\psi)$, set of test points ${\cal Z}$, containing $\bar x$,}\\
& \\
& \mbox{constant $L > 0$, and tolerance $\delta \geq 0$.}\\
& \\
\mbox{\bf Iteration:} & \mbox{Using the input data, form
the optimization problem (\ref{prob-ADual}) and }\\
& \\
& \mbox{compute its approximate solution $\bar
\lambda$ satisfying condition (\ref{eq-OptL}).}\\
& \\
\mbox{\bf Output:} & \mbox{Points $\bar s = L \nabla
d(\bar x) - G \bar \lambda$ and $x_+ = y^*_L(\bar s)$.}\\
& \\
\hline
\end{array}
\end{equation}
\begin{theorem}\label{th-Step}
1. Let point $x_+$ be generated by one iteration
(\ref{it-IGMM}) of the Inexact Gradient Method with Memory (IGMM),
and let $L \geq L_d(f)$. Then for any $y \in \dom
\psi$
\begin{equation}\label{eq-Prog}
\begin{array}{rcl}
\beta(x_+,y) & \leq & \beta(\bar x,y) + {1 \over L}
\left[ F(y) - F(x_+) + \delta \right];
\end{array}
\end{equation}
2. For every $y \in \dom \psi$ that satisfies $\beta_d(z_i,y) \geq \beta_d(\bar x, y)$, $i = 1,
\dots, m$ we have
\begin{equation}\label{eq-Prog1}
\begin{array}{rcl}
\beta(x_+,y) & \leq &
\left(1-{1 \over L}\, \mu_d(f)\right)\beta(\bar x,y) + {1 \over L}
\left[ F(y) - F(x_+) + \delta \right].
\end{array}
\end{equation}
\end{theorem}
\begin{proof}
Note that
$$
\begin{array}{rl}
& \beta(x_+,y) - \beta(\bar x, y) \\
\\
\refEQ{def-Breg} & d(y) - d(x_+) - \la \nabla d(x_+), y -
x_+ \ra - d(y) + d(\bar x) + \la \nabla d(\bar x), y -
\bar x \ra\\
\\
\refEQ{def-Breg} & \la \nabla d(\bar x) - \nabla d(x_+), y
- x_+ \ra - \beta(\bar x, x_+).
\end{array}
$$
The point $x_{+}$ is defined as
$$
x_+ = \arg\max\limits_{x \in \dom \psi} \left\{ \la
L \nabla d(\bar x) - G \bar \lambda , x \ra - L d(x) - \psi(x) \right\} .
$$
The first-order optimality condition for $x_+$ at point $y$ can be written in the following form:
$$
\begin{array}{rcl}
\psi(x_+) \leq & \psi(y) + \la G \bar \lambda + L (\nabla
d(x_+) - \nabla d(\bar x)), y - x_+ \ra.
\end{array}
$$
Hence,
$$
\begin{array}{rcl}
\beta(x_+,y) - \beta(\bar x, y) & \leq & {1 \over L}\left[
\psi(y) - \psi(x_+) + \la G \bar \lambda, y - x_+ \ra
\right] - \beta(\bar x, x_+).
\end{array}
$$
Note that
$$
\begin{array}{rcl}
\la G \bar \lambda, y - x_+ \ra & = & \la \bar \lambda,
G^T(y - x_+) \ra \; = \; \sum\limits_{i=1}^m \bar
\lambda^{(i)} \la g_i, y - x_+ \ra\\
\\
& = & \sum\limits_{i=1}^m \bar \lambda^{(i)} [ \la g_i,
z_i - x_+ \ra + \la g_i, y - z_i \ra ]\\
\\
& \refLE{def-RSmooth} & \sum\limits_{i=1}^m \bar
\lambda^{(i)} [ \la g_i, z_i - x_+ \ra + f(y) - f_i -
\mu_d(f) \beta_d(z_i,y)]\\
\\
& = & f(y) - \sum\limits_{i=1}^m
\bar \lambda^{(i)} [ f_i + \la g_i, x_+ - z_i\ra] - \mu_d(f) \sum\limits_{i=1}^m \bar
\lambda^{(i)}\beta_d(z_i,y).
\end{array}
$$
Under the conditions of Item 1, we drop the last term in the
above inequality and by Lemma~\ref{lm-Tol} obtain the
following:
$$
\begin{array}{rcl}
\beta(x_+,y) - \beta(\bar x, y) & \leq & {1 \over L}\left[
F(y) - \sum\limits_{i=1}^m \bar \lambda^{(i)} [ f_i + \la
g_i, x_+ - z_i\ra] - \psi(x_+) - L \beta(\bar x, x_+)
\right]\\
\\
& \leq & {1 \over L}\left[ F(y) - \max\limits_{1 \leq i
\leq m} [ f_i + \la g_i, x_+ - z_i\ra] - \psi(x_+) - L
\beta(\bar x, x_+) + \delta \right].
\end{array}
$$
Under the conditions of Item 2, by the same reasoning we get
$$
\begin{array}{l}
\beta(x_+,y) - \beta(\bar x, y)\\
\\
\leq \; {1 \over L}\left[ F(y) - \mu_d(f)\beta_d(\bar x,y)
- \sum\limits_{i=1}^m
\bar \lambda^{(i)} [ f_i + \la g_i, x_+ - z_i\ra] -
\psi(x_+) - L \beta(\bar x, x_+)
\right]\\
\\
\leq \; {1 \over L}\left[ F(y) - \mu_d(f)\beta_d(\bar x,y)
- \max\limits_{1 \leq i \leq m} [ f_i + \la g_i, x_+ -
z_i\ra] - \psi(x_+) - L \beta(\bar x, x_+) + \delta
\right].
\end{array}
$$
In both cases, since $\bar x \in {\cal Z}$, we have
$$
\begin{array}{rcl}
\max\limits_{1 \leq i \leq m} [ f_i + \la g_i, x_+ -
z_i\ra] + L \beta(\bar x, x_+) & \geq & f(\bar x) + \la
\nabla f(\bar x), x_+ - \bar x \ra + L \beta(\bar x, x_+)\\
\\
& \refGE{def-RSmooth} & f(x_+).
\end{array}
$$
Thus, we obtain inequalities (\ref{eq-Prog}) and
(\ref{eq-Prog1}).
\end{proof}
\begin{remark}\label{rm-Step}
In the end of the proof, we have seen that the statement
of Theorem~\ref{th-Step} remains valid if the
condition $L \geq L_d(f)$ is replaced by the following:
\begin{equation}\label{eq-Relax}
\begin{array}{rcl}
\max\limits_{1 \leq i \leq m} [ f_i + \la g_i, x_+ -
z_i\ra] + L \beta(\bar x, x_+) & \geq & f(x_+).
\end{array}
\end{equation}
\end{remark}

Denote the output of the iteration (\ref{it-IGMM}) by
$x_{\delta,L}(\bar x, {\cal Z})$. Then we can define the
following Inexact Gradient Method with Memory.
\begin{equation}\label{met-IGMM}
\begin{array}{|l|}
\hline \\
\mbox{{\bf Choose} $x_0 \in \inter (\dom \psi)$, $\delta \geq 0$ and $L>0$.}\\
\\
\mbox{\bf For $k \geq 0$, iterate:} \\
\\
\mbox{1). Choose the set ${\cal Z}_k$ containing $x_k$.}\\
\\
\mbox{2). Compute $x_{k+1} = x_{\delta,L}(x_k, {\cal
Z}_k)$.}\\
\\
\hline
\end{array}
\end{equation}

Let us describe the rate of convergence of this process.
\begin{theorem}\label{th-RateIGMMa}
Let sequence $\{x_k \}_{k \geq 1}$ be generated by IGMM
(\ref{met-IGMM}) with $L \geq L_d(f)$. Then, for any $T
\geq 1$ and $y \in \dom \psi$ we have
\begin{equation}\label{eq-RateA}
\begin{array}{rcl}
{1 \over T} \sum\limits_{k=1}^T F(x_k) & \leq & F(y) + {L
\over T} \beta_d(x_0,y) + \delta.
\end{array}
\end{equation}
\end{theorem}
\begin{proof}
Indeed, in view of inequality (\ref{eq-Prog}), we have
$$
\begin{array}{rcl}
\beta_d(x_{k+1},y) & \leq & \beta_d(x_k,y) + {1 \over L}
\left[ F(y) - F(x_{k+1}) + \delta \right], \quad k \geq 0.
\end{array}
$$
Summing up these inequalities for $k = 0, \dots, T-1$, we
get inequality (\ref{eq-RateA}).
\end{proof}

In the above result, the only restriction for the sets
${\cal Z}_k$ is the inclusion (\ref{eq-Inc}). If we apply
a more accurate strategy of choosing ${\cal Z}_k$, we can
get for this scheme a finer estimate of its rate of
convergence.
\begin{theorem}\label{th-RateIGMMb}
Let sequence $\{x_k \}_{k \geq 1}$ be generated by IGMM
(\ref{met-IGMM}) with $L \geq L_d(f)$. Assume that besides
the condition (\ref{eq-Inc}), the sets ${\cal Z}_k$
satisfy also the following condition:
\begin{equation}\label{eq-CondZ}
\begin{array}{rcl}
{\cal Z}_k & \subseteq & \{x_0, \dots, x_k\}, \quad k \geq
0.
\end{array}
\end{equation}
Suppose that $\mu_d(f)>0$ and for all $k$, $1 \leq k \leq
T$, we have
\begin{equation}\label{eq-BigF}
F(x_k) - F^* \geq \delta.
\end{equation}
Then for $\Delta^*_T = \min\limits_{1 \leq k \leq T}
\left[F(x_k) - F^*\right]$ we get the following rate of
convergence
\begin{equation}\label{eq-RateB}
\begin{array}{rcl}
\Delta^*_T & \leq & \delta + { (1-\gamma)^T \mu_d(f) \over
1 - (1-\gamma)^T}\, \beta_d(x_0,x^*) \\
\\
& \leq & \delta + { \mu_d(f) \over e^{\gamma T} - 1}\,
\beta_d(x_0,x^*) \; \leq \; \delta + {L \over T}
\beta_d(x_0,x^*),
\end{array}
\end{equation}
where $\gamma = {1 \over L} \, \mu_d(f)$.
\end{theorem}
\begin{proof}
In view of assumption (\ref{eq-BigF}) and inequality
(\ref{eq-Prog}), we have
$$
\begin{array}{rcl}
\beta_d(x_{k+1},x^*) & \leq & \beta_d(x_k,x^*), \quad 0
\leq k \leq T-1.
\end{array}
$$
Therefore, in view of inequality (\ref{eq-Prog1}), for
$r_k \Def \beta_d(x_k,x^*)$ and all $k = 0, \dots, T-1$ we
have
$$
\begin{array}{rcl}
r_{k+1} & \leq & (1-\gamma)r_k - {1 \over L} [F(x_{k+1}) -
F^* - \delta] \; \leq \; (1-\gamma)r_k - {1 \over L}
[\Delta^*_T - \delta].
\end{array}
$$
Applying this inequality recursively, we get
$$
\begin{array}{rcl}
{1 \over L} [\Delta^*_T - \delta] {1 - (1-\gamma)^T \over
1 - (1-\gamma)} & \leq & (1-\gamma)^T r_0.
\end{array}
$$
This can be rewritten as
$$
\begin{array}{rcl}
\Delta^*_T & \leq & \delta + {L \gamma (1-\gamma)^T \over
1 - (1-\gamma)^T} \beta_d(x_0,x^*) \; = \; \delta + {
(1-\gamma)^T \mu_d(f) \over 1 - (1-\gamma)^T}
\beta_d(x_0,x^*). \hspace{5ex}
\end{array}
$$
\end{proof}

Note that the rate of convergence given by inequality
(\ref{eq-RateB}) is continuous as $\mu_d(f) \to 0$.

As we have mentioned in Remark \ref{rm-L}, it is important
to adjust the value of the constant $L$ during the
minimization process. Therefore we present an adaptive
version of the method (\ref{met-IGMM}).
\begin{equation}\label{met-IGMMA}
\begin{array}{|l|}
\hline \\
\mbox{{\bf Choose} $x_0 \in \inter (\dom \psi)$, $\delta
\geq 0$, and some
$L_0 \in (0, L_d(f)]$.}\\
\\
\mbox{\bf For $k \geq 0$, iterate:} \\
\\
\mbox{1). Choose the set ${\cal Z}_k$ containing $x_k$.}\\
\\
\mbox{2). Find the smallest integer $i_k \geq 0$ such
that}\\
\\
\mbox{for the point $x^+_{k} = x_{\delta,2^{i_k}L_k}(x_k,
{\cal
Z}_k)$ we have}\\
\\
\hspace{10ex} f(x^+_k) \refLE{eq-Mod} \ell_k(x^+_k) + L_k
\beta_d(x_k,
x^+_k).\\
\\
\mbox{3). Set $x_{k+1} = x^+_k$ and $L_{k+1} =
2^{i_k-1}L_k$.}\\
\\
\hline
\end{array}
\end{equation}

The rate of convergence of this algorithm can be
established exactly in the same way as the one of method
(\ref{met-IGMM}). The main fact is that during the
minimization process we always have
$$
L_k \leq 2 L_d(f), \quad k \geq 0.
$$
Therefore, for any $y \in \dom \psi$ we will have
$$
\begin{array}{rcl}
\beta_d(x_{k+1},y) & \leq & \beta_d(x_k,y) + {1 \over 2
L_d(f)} [F(y) - F(x_{k+1}) + \delta]
\end{array}
$$
with corresponding consequences for the rate of
convergence. At the same time, the average number oracle
calls at each iteration of this method is bounded by
two (see~\cite{ref_007} for justification details).

\section{Getting an approximate solution of the anti-dual
problem}\label{sc-Approx}

The complexity of solving the auxiliary problem
(\ref{prob-ADual}) crucially depends on the properties of
the prox-function $d(\cdot)$. In the previous section, we
assumed its strict convexity and the solvability of problem
(\ref{def-Phi}) (see Assumption \ref{ass-Phi}). It is time
now to make a stronger assumption, which ensures these two
properties.
\begin{assumption}\label{ass-Strong}
Function $d(\cdot)$ is differentiable in the interior of
its domain and strongly convex with convexity parameter
one:
\begin{equation}\label{eq-Strong}
\begin{array}{rcl}
d(y) & \geq & d(x) + \la \nabla d(x), y - x \ra + \half \|
y - x \|^2, \quad x \in \inter (\dom d), \; y \in \dom d.
\end{array}
\end{equation}
\end{assumption}
Clearly, for all $x, y \in \inter (\dom d)$ we have
\begin{equation}\label{eq-BStrong}
\begin{array}{rcl}
\beta_d(x,y) & \geq & \half \| y - x \|^2.
\end{array}
\end{equation}

The main consequence of Assumption \ref{ass-Strong} is
the Lipschitz continuity of the gradient of function
$\Phi(\cdot)$. Since usually this fact is proved for a
function $\psi(\cdot)$ being an indicator function of a
closed convex set, we provide it with a simple proof.
\begin{lemma}\label{lm-Lip}
Let function $d(\cdot)$ satisfy Assumption
\ref{ass-Strong}. Then the gradient $\nabla \Phi_L(s) =
y^*_L(s)$, $s \in \E^*$, is Lipschitz continuous:
\begin{equation}\label{eq-LipGP}
\begin{array}{rcl}
\| \nabla \Phi_L(s_1) - \nabla \Phi_L(s_2) \| & \leq & {1
\over L} \| s_1 - s_2 \|_*, \quad s_1, s_2 \in \E^*.
\end{array}
\end{equation}
\end{lemma}
\proof
Let us write down the first-order optimality conditions
for the optimization problems defining the points $y_1 \Def
y^*_L(s_1)$ and $y_2 \Def y^*_L(s_2)$:
$$
\begin{array}{rcl}
\la s_1 - L[\nabla d(y_1) - \nabla d(\bar x)], y - y_1 \ra
- \psi(y) & \leq - \psi(y_1), \quad y \in \dom
\psi,\\
\\
\la s_2 - L[\nabla d(y_2) - \nabla d(\bar x)], y - y_2 \ra
- \psi(y) & \leq - \psi(y_2), \quad y \in \dom
\psi.
\end{array}
$$
Taking in the first inequality $y = y_2$ and $y = y_1$ in
the second one, and adding the results, we obtain
$$
\begin{array}{rcl}
\la s_1 - s_2, y_2 - y_1 \ra & \leq & L \la \nabla d(y_1)
- \nabla d(y_2), y_2 - y_1 \ra.
\end{array}
$$
Thus,
$$
\begin{array}{rcl}
\la s_1 - s_2, y_1 - y_2 \ra & \geq & L \la \nabla d(y_2)
- \nabla d(y_1), y_2 - y_1 \ra \; \refGE{eq-Strong} L \|
y_1 - y_2 \|^2.
\end{array}
$$
Therefore, by the Cauchy-Schwartz inequality, we get
$$
\begin{array}{rcl}
\| y^*_L(s_1) - y^*_L(s_2) \| & \leq & {1 \over L} \| s_1
- s_2 \|_*. \hspace{5ex} \Box
\end{array}
$$

Thus, in this section, our main problem of interest is as
follows:
\begin{equation}\label{prob-ADual1}
\xi^*_L \; = \; \min\limits_{\lambda \in \Delta_m}
\left\{ \xi_L(\lambda) \Def \Phi_L(L\nabla d(\bar x) - G
\lambda) + \la \lambda, f_* \ra + \alpha \right\},
\end{equation}
where $\alpha = L [ d(\bar x) - \la \nabla d(\bar x), \bar
x \ra]$. This is a convex optimization problem over a
simplex, where the objective function has a Lipschitz-continuous
gradient.

The most natural algorithm for solving the problem
(\ref{prob-ADual}) is the {\em Frank-Wolfe algorithm}
\cite{ref_003} (or {\em Conditional Gradients Method}). For our
problem, it looks as follows.
\begin{equation}\label{met-FW}
\begin{array}{|l|}
\hline \\
\mbox{{\bf Set} $\lambda_0 = {1 \over m} \bar e_m$.}\\
\\
\mbox{\bf For $k \geq 0$ iterate:}\\
\\
\mbox{1. Compute the gradient $\nabla \xi_L(\lambda_k)$.}\\
\\
\mbox{2. Compute $i_k = \arg\min\limits_{1 \leq i \leq m}
\nabla_i \xi_L(\lambda_k)$.}\\
\\
\mbox{3. Set $\lambda_{k+1} = {k \over k+2} \lambda_k + {2
\over k+2} e_{i_k}$.}\\
\\
\hline
\end{array}
\end{equation}
In this scheme, $\bar e_m \in \R^m$ is the vector of all
ones, and $e_i$ is $i$th coordinate vector in $\R^m$.

In order to estimate the rate of convergence of this
method, we introduce the following accuracy measure:
$$
\begin{array}{rcl}
\delta_L(\bar \lambda) & = & \max\limits_{\lambda \in
\Delta_m} \la \nabla \xi_L(\bar \lambda), \bar \lambda -
\lambda \ra.
\end{array}
$$
For the sequence $\{ \lambda_k \}_{k \geq 0}$ generated by
the method (\ref{met-FW}), denote
$$
\begin{array}{rcl}
\delta^*_L(T) & = & \min\limits_{0 \leq k \leq T}
\delta_L(\lambda_k), \quad T \geq 0.
\end{array}
$$
For estimating the rate of convergence of method
(\ref{met-FW}), we need to choose an appropriate norm in
$\R^m$. Since the feasible set of the problem
(\ref{prob-ADual}) is the standard simplex, it is
reasonable to use the $\ell_1$-norm:
$$
\begin{array}{rcl}
\| \lambda \|_1 & = & \sum\limits_{i=1}^m |
\lambda^{(i)}|, \quad \lambda \in \R^m.
\end{array}
$$
Then, for measuring the gradients of function $\xi_L(\cdot)$,
we can use the $\ell_{\infty}$-norm:
$$
\begin{array}{rcl}
\| \lambda \|_{\infty} & = & \max\limits_{1 \leq i \leq m}
| \lambda^{(i)}|, \quad \lambda \in \R^m.
\end{array}
$$
In this case, the Lipschitz constant for the gradients of
function $\xi_L(\cdot)$ can be estimated as follows:
$$
\begin{array}{cl}
& \| \nabla \xi_L(\lambda_1) - \nabla \xi_L(\lambda_2)
\|_{\infty} \\
\\
= & \max\limits_{1 \leq i \leq m} | \la g_i , \nabla
\Phi(L \nabla d(\bar x) - G \lambda_2) -
\nabla \Phi(L \nabla d(\bar x) - G \lambda_1) \ra| \\
\\
\leq & \max\limits_{1 \leq i \leq m} \| g_i \|_* \cdot \|
\nabla \Phi(L \nabla d(\bar x) - G \lambda_2) - \nabla
\Phi(L \nabla d(\bar x) - G \lambda_1) \|\\
\\
\refLE{eq-LipGP} & \max\limits_{1 \leq i \leq m} \| g_i
\|_* \cdot {1 \over L} \| G(\lambda_1-\lambda_2) \|_* \;
\leq \; {1 \over L} \max\limits_{1 \leq i \leq m} \| g_i
\|^2_* \cdot \| \lambda_1 - \lambda_2 \|_1.
\end{array}
$$
Thus, the gradients of function $\xi_L(\cdot)$ are
Lipschitz continuous with the constant
\begin{equation}\label{def-LipGX}
\begin{array}{rcl}
L(\xi_L) & = & {1 \over L} \max\limits_{1 \leq i \leq m}
\| g_i \|^2_*.
\end{array}
\end{equation}
Since the diameter of the standard simplex in $\R^m$ in
$\ell_1$-norm is two, in accordance to the estimate (3.13)
in \cite{ref_006}, we can guarantee the following rate of
convergence:
\begin{equation}\label{eq-RateFW}
\begin{array}{rcl}
\delta^*_L(T) & \leq & {18 \over L \cdot T}
\max\limits_{1 \leq i \leq m} \| g_i \|^2_* ,
\quad t \geq 1.
\end{array}
\end{equation}
(Here we replace the constant ${136 \over 11 \ln 2}$ from
\cite{ref_006} by a bigger value $18$.) In accordance to the
condition (\ref{eq-OptL}), this means that we need
\begin{equation}\label{eq-ItFW}
\begin{array}{rcl}
N_L(\delta) & = & {18 \over L \cdot \delta} \max\limits_{1
\leq i \leq m} \| g_i \|^2_*
\end{array}
\end{equation}
iterations of the method (\ref{met-FW}) in order to
generate an appropriate dual solution $\bar \lambda$.

\section{Unconstrained minimization in Euclidean
setup}\label{sc-Euclid}

In this section we consider the simplest unconstrained
minimization problem
\begin{equation}\label{prob-Un}
\begin{array}{rcl}
f^* & = & \min\limits_{x \in \E} f(x),
\end{array}
\end{equation}
where $f(\cdot)$ is a smooth convex function. For
measuring distances in $\E$, we introduce a Euclidean norm
$$
\begin{array}{rcl}
\| x \| & = & \la B x, x \ra^{1/2}, \quad x \in \E,
\end{array}
$$
where $B = B^* \succ 0$ is a linear operator from $\E$ to
$\E^*$. Then the dual norm is defined as follows:
$$
\begin{array}{rcl}
\| g \|_* & = & \la g, B^{-1} g \ra^{1/2}, \quad g \in
\E^*.
\end{array}
$$
Let us choose now the distance function $d(x) = \half \| x
\|^2$. Then the Bregman distance is given by
$$
\begin{array}{rcl}
\beta_d(x,y) & = & \half \| x - y \|^2, \quad x, y \in \E.
\end{array}
$$
In this case, the relative smoothness condition
(\ref{def-RSmooth}) is equivalent to strong convexity and
Lipschitz continuity of the gradient:
\begin{equation}\label{eq-LS}
\begin{array}{rcll}
f(y) - f(x) - \la \nabla f(x), y - x \ra & \geq & \half
\mu_d(f)
\| x - y \|^2,\\
& & & \quad x,y \in \dom f,\\
f(y) - f(x) - \la \nabla f(x), y - x \ra & \leq & \half
L_d(f) \| x - y \|^2.
\end{array}
\end{equation}

Let us write down now the specific form of the objective
function $\xi_L(\cdot)$ in problem (\ref{prob-ADual}).
Note that in our case
$$
\begin{array}{rcl}
\Phi_L(s) & = & \max\limits_{y \in \E} \left\{ \la s, y \ra -
{L \over 2} \| y \|^2 \right\} \; = \; {1 \over 2L} \| s \|_*^2, \quad
s \in \E^*.
\end{array}
$$
Therefore,
$$
\begin{array}{rcl}
\xi_L(\lambda) = {1 \over 2L} \| L B \bar x - G \lambda
\|^2_* + \la \lambda, f_* \ra + \alpha,
\end{array}
$$
where $\alpha = - \half L \| \bar x \|^2$. The gradient of
function $\xi_L(\cdot)$ can be computed as follows:
\begin{equation}\label{eq-Grad}
\begin{array}{rcl}
\nabla \xi_L(\lambda) & = & {1 \over L} G^* B^{-1} (G
\lambda - L B
\bar x) + f_*, \quad \lambda \in \R^m\\
\\
& = & {1 \over L} Q \lambda - \bar f,
\end{array}
\end{equation}
where $Q = G^* B^{-1} G$ and $\bar f = G^* \bar x - f_*$.
Note that
$$
\begin{array}{rcl}
\bar f^{(i)} & \refEQ{def-f*} & f_i + \la g_i, \bar x - z_i \ra, \quad i =
1, \dots, m.
\end{array}
$$
Thus, in the Euclidean setup, our auxiliary problem
(\ref{prob-ADual}) can be written as follows:
\begin{equation}\label{prob-EDual}
\begin{array}{rcl}
\min\limits_{\lambda \in \Delta_m} \Big[ \; \xi_L(\lambda)
& = & {1 \over 2L} \la \lambda, Q \lambda \ra - \la
\lambda,
\bar f \ra \; \Big].
\end{array}
\end{equation}
The stopping criterion (\ref{eq-OptL}) for this problem is
as follows:
$$
\begin{array}{rcl}
\la \bar \lambda, \nabla \xi_L(\bar \lambda) \ra &
\refEQ{eq-Grad} & \la \bar \lambda, {1 \over L} Q \bar
\lambda - \bar f \ra \; \refLE{eq-OptL} \; \delta +
\min\limits_{1 \leq i \leq m}
\nabla_i \xi_L(\bar \lambda)\\
\\
& = & \delta + \min\limits_{1 \leq i \leq m} \left({1 \over L}
Q \bar \lambda - \bar f\right)^{(i)}.
\end{array}
$$
Note that the main output of the minimization process for
problem (\ref{prob-EDual}) is
$$
\begin{array}{rcl}
x_+ & = & y_*(L B \bar x - G \bar \lambda ) \;
\refEQ{eq-GPhi} \; {1 \over L} B^{-1} (L B \bar x - G \bar
\lambda ) \; = \; \bar x - {1 \over L} B^{-1} G \bar
\lambda.
\end{array}
$$
Then
$$
\begin{array}{rcl}
\bar f - {1 \over L} Q \bar \lambda & = & \bar f - {1
\over L} G^* B^{-1} G \bar \lambda \; = \; \bar f +
G^*(x_+ - \bar x).
\end{array}
$$
Hence, in the Euclidean case, the stopping criterion
(\ref{eq-OptL}) can be written as follows:
\begin{equation}\label{eq-AStop}
\begin{array}{rcl}
\sum\limits_{i=1}^m \bar \lambda^{(i)} [ f_i + \la g_i,
x_+ - z_i \ra] & \leq & \delta + \max\limits_{1 \leq i
\leq m} [ f_i + \la g_i, x_+ - z_i \ra].
\end{array}
\end{equation}

Now we can estimate the computational expenses of the
method (\ref{met-FW}) as applied to the auxiliary problem
(\ref{prob-EDual}).
\begin{enumerate}
\item
Computation of the matrix $Q$: $O(m^2 n)$ arithmetic
operations. For certain strategies for {\em updating} the
sets ${\cal Z}_k$, it can be reduced to $O(mn)$
operations.
\item
Computation of the vector $\bar f$: $O(mn)$ operations.
\item
Computation of the initial gradient $u_0 = {1 \over L}
\lambda_0 - \bar f$: $O(m^2)$ operations. For certain
updating strategies it can be $O(m)$.
\item
Expenses at each iteration:
\begin{itemize}
\item
Computing the index $i_k$: $O(m)$ operations.
\item
Updating the point $\lambda_k$: $O(m)$ operations.
\item
Updating the gradient $u_k = {1 \over L} Q \lambda_k -
\bar f$: $O(m)$ operations.
\end{itemize}
\end{enumerate}

Thus, taking into account the upper bound
(\ref{eq-ItFW}) for the number of iterations in method
(\ref{met-FW}), we obtain the following bound for the
arithmetic complexity of problem (\ref{prob-EDual}) with
reasonable updating strategies for the sets ${\cal Z}_k$:
\begin{equation}\label{eq-EComp}
\begin{array}{c}
O \left( m n + {m \over L \cdot \delta} \max\limits_{1
\leq i \leq m} \| g_i \|^2_* \right).
\end{array}
\end{equation}
Taking into account that we can expect that in the problem
(\ref{prob-Un}) we have
$$
\begin{array}{rcl}
{1 \over 2L_d(f)} \| \nabla f(z_i) \|^2_* & \refLE{eq-LS} &
f(z_i) - f^* \; \to \; 0
\end{array}
$$
as $i \to \infty$, the bound in (\ref{eq-EComp}) suggests that
the overhead of solving the inner problem (\ref{prob-EDual}) decreases
to a small constant in $O(m n)$ as the algorithm approaches the optimum.

\section{Numerical experiments}\label{sc-Num}

In this section we present preliminary computational
results for method (\ref{met-IGMMA}) as applied to the
following unconstrained minimization problem:
\begin{equation}\label{prob-UTest}
\begin{array}{rcl}
\min\limits_{x \in \R^n} \Big[ \; f(x) & = & \mu \ln
\left( \sum\limits_{j=1}^M e^{(\la a_j , x \ra - b_j)/\mu}
\right) \Big].
\end{array}
\end{equation}
The data defining this function is randomly generated in
the following way. First of all, we generate a collection
of random vectors
$$
\hat a_1, \dots, \hat a_M
$$
with entries uniformly distributed in the interval
$[-1,1]$. Using the same distribution, we generate values
$b_j$, $j = 1, \dots, m$. Using this data, we form the
preliminary function
$$
\begin{array}{rcl}
\hat f(x) & = & \mu \ln
\left( \sum\limits_{j=1}^M e^{(\la a_j , x \ra - b_j)/\mu}
\right)
\end{array}
$$
and compute $g = \nabla \hat f(0)$. Then, we define
$$
\begin{array}{rcl}
a_j = \hat a_j - g, \quad j = 1, \dots, n.
\end{array}
$$
Clearly, in this case we have $\nabla f(0) = 0$, so the
unique solution of our test problem (\ref{prob-UTest}) is
$x^* = 0$. The starting point $x_0$ is chosen in
accordance to the uniform distribution on the Euclidean
sphere of radius one.

Thus, the problem (\ref{prob-UTest}) has three parameters,
the dimension $n$, the number of linear functions $M \geq n$,
and the smoothness coefficient $\mu>0$. In our
experiments, we always choose $M = 6n$. Let us present our
computational results for different values of $n$ and
$\mu$.

In the definition of method (\ref{met-IGMMA}) we have some
freedom in the choice of the bundle ${\cal Z}_k$. Let us
bound its maximal size by a parameter $m \geq 1$. Then $m
= 1$ corresponds to the usual Gradient Method. In the
first series of our experiments (shown in Tables \ref{label_001}
and \ref{label_002}) we always choose
$$
m = n.
$$

We also have some freedom in the updating strategy for the
sets ${\cal Z}_k$. Clearly, at the first $m$ steps we can
simply add all new points in the bundle. However, at the
next iterations we need to decide on the strategy of
replacement of the the old information. In our experiments
we implemented two strategies:
\begin{itemize}
\item
Cyclic replacement ({\sc Cyclic}).
\item
Replacement of the linear function with the maximal norm of the
gradient ({\sc Max-Norm}).
\end{itemize}
The second strategy is motivated by the formula
(\ref{def-LipGX}) for the Lipschitz constant of the
gradient of function $\xi_L(\cdot)$. For both strategies,
at each iteration we need to update only one column of
matrix $Q_k$ (see (\ref{eq-Grad})), which costs $O(mn)$
operations.

Let us present the results of our numerical experiments.
All methods were stopped when the residual in the function
value was smaller than $\epsilon = 10^{-6}$. The parameter
$\delta$ for the stopping criterion (\ref{eq-OptL}) was
chosen as $\delta = \epsilon/2$.

In Tables \ref{label_001} and \ref{label_002}, the first line
indicates the total number of iterations. The second line
displays the total number of oracle calls.
The third line shows the average
number of Frank-Wolfe steps per iterations (for the Gradient
Method we just put two). The next line indicates the total
computational time (in seconds). Finally, at the last line
we can see the average time spent on one iteration of the
corresponding method (in milliseconds).

\begin{table}
\caption{Smoothness parameter $\mu = 0.05$.}
\centering
\small
{\begin{tabular}{|l|r|r|r|}
\hline
$n = 100$ & GM & Cyclic & Max-Norm \\
\hline
Iter & 2683 & 801 & 664 \\
NFunc & 5371 & 1606 & 1332 \\
FW/Iter & 2 & 55 & 66 \\
Time (s) & 1.94 & 0.88 & 0.73 \\
IT(ms) & 0.72 & 1.10 & 1.10 \\
\hline
\end{tabular}}
{\begin{tabular}{|l|r|r|r|}
\hline
$n = 250$ & GM & Cyclic & Max-Norm \\
\hline
Iter & 2148 & 227 & 227 \\
NFunc & 4302 & 459 & 459 \\
FW/Iter & 2 & 243 & 243 \\
Time (s) & 10.20 & 1.36 & 1.36 \\
IT(ms) & 4.75 & 5.99 & 5.99 \\
\hline
\end{tabular}}
\\[2mm]
{\begin{tabular}{|l|r|r|r|}
\hline
$n = 500$ & GM & Cyclic & Max-Norm \\
\hline
Iter & 2902 & 268 & 268 \\
NFunc & 5809 & 537 & 537 \\
FW/Iter & 2 & 428 & 428 \\
Time (s) & 52.16 & 5.5 & 5.5 \\
IT(ms) & 17.97 & 20.52 & 20.52 \\
\hline
\end{tabular}}
\label{label_001}
\end{table}

\begin{table}
\caption{Smoothness parameter $\mu = 0.01$}
\centering
\small
{\begin{tabular}{|l|r|r|r|}
\hline
$n = 100$ & GM & Cyclic & Max-Norm \\
\hline
Iter & 43893 & 4171 & 6710 \\
NFunc & 87795 & 8351 & 13427 \\
FW/Iter & 2 & 59 & 17 \\
Time (s) & 31.73 & 4.50 & 6.95 \\
IT(ms) & 0.72 & 1.08 & 1.04 \\
\hline
\end{tabular}}
{\begin{tabular}{|l|r|r|r|}
\hline
$n = 250$ & GM & Cyclic & Max-Norm \\
\hline
Iter & 116479 & 45183 & 25492 \\
NFunc & 232967 & 90377 & 50990 \\
FW/Iter & 2 & 8 & 13 \\
Time (s) & 540.64 & 311.95 & 176.02 \\
IT(ms) & 4.64 & 6.90 & 6.90 \\
\hline
\end{tabular}}
\\[2mm]
{\begin{tabular}{|l|r|r|r|}
\hline
$n = 500$ & GM & Cyclic & Max-Norm \\
\hline
Iter & 105610 & 38144 & 29916 \\
NFunc & 211229 & 76297 & 59840 \\
FW/Iter & 2 & 10 & 17 \\
Time (s) & 1894.19 & 1010.30 & 788.97 \\
IT(ms) & 17.94 & 26.49 & 26.37 \\
\hline
\end{tabular}}
\label{label_002}
\end{table}

As we can see from these tables, in all our experiments
the gradient methods with memory were better than the
standard Gradient Method, both in the number of
iterations, and, what is rather surprising, in the total
computational time. The {\sc Max-Norm} version usually
outperforms the {\sc Cyclic} version. It is interesting
that the auxiliary algorithm~(\ref{met-FW}) works very
well. The average time spent on one iteration of the
methods with memory is never increased more than on $50\%$
of the time of the simple Gradient Method. This is
partially explained by the fact that in our test problems
the data is fully dense, so each call of oracle is very
expensive ($O(Mn)$ operations).

Let us look now at how small bundles can accelerate the
Gradient Method. In Tables \ref{label_003} and \ref{label_004},
the first line with parameter Bundle = 1 corresponds to the
Gradient method with line search. The next lines display the results
for different sizes of the bundle. We list the number of
iterations, the average number of Frank-Wolfe steps per
iteration, and the total computational time in seconds.
Table~\ref{label_003} displays the results for the IGMM
(\ref{met-IGMMA}) with the cyclic replacement strategy for each
bundle size. In Table~\ref{label_004}, we show the results for the Max-Norm
replacement strategy. The accuracy parameters for the experiments shown
in Tables \ref{label_003} and \ref{label_004} are $\epsilon = 10^{-4}$ and
$\delta = \epsilon/2$. The smoothness parameter for our
objective function is chosen as $\mu = 0.05$.

\begin{table}
\caption{Gradient Method with Cyclic Memory Replacement}
\centering
\small
{\begin{tabular}{|r|rrr|rrr|rrr|} \hline
& \multicolumn{3}{c|}{$n = 100$} &
\multicolumn{3}{c|}{$n = 200$} & \multicolumn{3}{c|}{$n = 400$} \\ \hline
\multicolumn{1}{|c|}{Bundle} & Iter & FW & Sec &
Iter & FW & Sec & Iter & FW & Sec \\ \hline
1 & 2683 & 2 & 1.94 & 1753 & 2 & 4.78 & 1676 & 2 & 19.86 \\
2 & 2543 & 3 & 1.84 & 1433 & 7 & 3.91 & 1669 & 5 & 19.81 \\
4 & 1755 & 13 & 1.31 & 829 & 39 & 2.31 & 1084 & 18 & 12.92 \\
8 & 1363 & 19 & 1.03 & 633 & 60 & 1.78 & 758 & 37 & 9.08 \\
16 & 1220 & 22 & 0.97 & 579 & 68 & 1.67 & 720 & 44 & 8.69 \\
32 & 1202 & 24 & 1.00 & 593 & 71 & 1.76 & 719 & 47 & 8.78 \\
64 & 1138 & 28 & 1.06 & 423 & 114 & 1.39 & 636 & 58 & 8.02 \\
128 & 880 & 60 & 1.06 & 247 & 320 & 0.94 & 349 & 127 & 4.63 \\
256 & 207 & 361 & 0.36 & 173 & 507 & 0.70 & 202 & 260 & 2.70 \\ \hline
\end{tabular}}
\label{label_003}
\end{table}

\begin{table}
\caption{Gradient Method with Max-Norm Memory Replacement}
\centering
\small
{\begin{tabular}{|r|rrr|rrr|rrr|} \hline
& \multicolumn{3}{c|}{$n = 100$} &
\multicolumn{3}{c|}{$n = 200$} & \multicolumn{3}{c|}{$n = 400$} \\ \hline
\multicolumn{1}{|c|}{Bundle} & Iter & FW & Sec &
Iter & FW & Sec & Iter & FW & Sec \\ \hline
1 & 2683 & 2 & 1.94 & 1753 & 2 & 4.76 & 1676 & 2 & 19.86 \\
2 & 1288 & 4 & 0.95 & 737 & 11 & 2.05 & 796 & 7 & 9.47 \\
4 & 592 & 27 & 0.45 & 316 & 79 & 0.91 & 385 & 38 & 4.61 \\
8 & 400 & 48 & 0.33 & 269 & 121 & 0.78 & 283 & 86 & 3.42 \\
16 & 362 & 58 & 0.33 & 418 & 91 & 1.22 & 335 & 86 & 4.05 \\
32 & 786 & 33 & 0.67 & 354 & 116 & 1.08 & 354 & 89 & 4.34 \\
64 & 838 & 39 & 0.80 & 379 & 127 & 1.23 & 421 & 85 & 5.31 \\
128 & 529 & 100 & 0.66 & 285 & 277 & 1.09 & 352 & 127 & 4.64 \\
256 & 207 & 361 & 0.36 & 173 & 507 & 0.70 & 202 & 260 & 2.70 \\ \hline
\end{tabular}}
\label{label_004}
\end{table}

As we can see from these tables, the Max-Norm replacement
strategy was always better than the cyclic one. Even for small
bundle sizes, the total number of iterations decreases very quickly.
What is more important, this decrease is also seen in total computation
time. The number of auxiliary Frank-Wolfe steps remains on an acceptable
level and cannot increase significantly the computational
time of each iteration as compared with the Gradient
Method. Recall that our test function has an expensive
oracle, requiring $O(Mn)$ operations for computing the
function value and the gradient. For the Max-Norm version of
IGMM, the optimal size of the bundle is probably between 8
and 16. Another candidate is 256, but it needs many more
Frank-Wolfe steps.

Maybe our preliminary conclusions are problem specific.
However, we believe that in any case they demonstrate a
high potential of our approach in increasing the efficiency of
gradient methods, both in accelerated and, hopefully,
non-accelerated variants.

\section*{Funding}

This project has received funding from the European Research Council (ERC) under the European Union's Horizon 2020 research and innovation programme (grant agreement No. 788368).


\begin{thebibliography}{NN}

\bibitem{ref_001}
H. Bauschke, J. Bolte, and M. Teboulle,
{\em A descent lemma beyond Lipschitz gradient continuity: First-order methods revisited and applications},
Math. Oper. Res. 42 (2017), pp.~330--348.

\bibitem{ref_002}
D. Drusvyatskiy, M. Fazel, and S. Roy,
{\em An optimal first order method based on optimal quadratic averaging},
SIAM J. Optim. 28 (2018), pp.~251--271.

\bibitem{ref_003}
M. Frank and P. Wolfe,
{\em An algorithm for quadratic programming},
Naval Res. Logist. Q. 3 (1956), pp.~149--154.

\bibitem{ref_004}
H. Lu, R. Freund, and Yu. Nesterov,
{\em Relatively smooth convex optimization by first-order methods, and Applications},
SIAM J. Optim. 28S (2018), pp.~333--354.

\bibitem{ref_005}
Yu. Nesterov,
{\em Gradient methods for minimizing composite functions},
Math. Program. 140 (2013), pp.~125--161.

\bibitem{ref_006}
Yu. Nesterov,
{\em Complexity bounds for primal-dual methods minimizing the model of objective function},
Math. Program. 171 (2018), pp.~311--330.

\bibitem{ref_007}
Yu. Nesterov,
{\em Lectures on Convex Optimization},
Springer, Berlin, Germany, 2018.

\end{thebibliography}
\end{document}